# A Note on the Stability of Exponential Dichotomy of Linear Differential Equations

**Osvaldo Mendez, Nada al Hanna**
**Department of Mathematical Sciences, University of Texas – El Paso, USA**
**124 Bell Hall, 500W University Ave., El Paso, TX, 79968 USA**
**e-mail:** mendez@math.utep.edu, nfalhanna@miners.utep.edu

**ABSTRACT.** We present an elementary Functional Analytic proof of the roughness of Exponential Dichotomy of Ordinary Differential Equations (with exponential growth) on an arbitrary Banach Space.
**KEYWORDS:** Exponential Dichotomy, Exponential Growth, Ordinary Differential Equation, Unbounded Operator, Banach Space, Roughness, Perturbation.

## 1 Introduction

We consider a Banach Space $X$ and a strongly continuous function:

$$A: \Re \to B(X) = \{T : X \to X, T \text{ linear and bounded}\}$$

and the non-autonomous differential equation

$$\frac{d}{dt}x = A(t)x \tag{1}$$

We denote the Cauchy operator of equation (1) by $U(t)$; in the sequel we assume that (1) has exponential growth, that is, that there exist positive constants $\alpha$ and $\beta$ (which will be referred to as the exponential growth constants) such that for any $t_0, t \in \Re$, one has

$$\left\| U(t)U^{-1}(t_0) \right\| \le \alpha\, e^{\beta(t-t_0)} \tag{2}$$

Equation (1) is said to have exponential dichotomy if there exist two bounded projections $P : X \to X$ and $Q : X \to X$ with $P + Q = I$ and positive constants $N_i$, $v_i$, $i = 1, 2$ for which the following estimates hold for any $t, s \in I$:

$$\left\|U(t)PU^{-1}(s)\right\| \leq N_1 \, e^{-v_1(t-s)} \text{ if } t \geq s \qquad (3)$$

$$\left\|U(t)QU^{-1}(s)\right\| \leq N_2 \, e^{-v_2(s-t)} \text{ if } s \geq t \qquad (4)$$

Equivalently, (1) is Exponentially Dichotomic if the subspace of $X$ consisting of initial values of solutions that are bounded on $[0,\infty)$ is complemented in $X$. It follows immediately from the definition that if (1) possesses Exponential Dichotomy, then any solution $x(t)$ that is bounded on $[0,\infty)$ tends, satisfies

$$\lim_{t \to \infty} x(t) = 0$$

and

$$\lim_{t \to -\infty} x(t) = \infty \, ;$$

analogously a solution bounded on the left semi-axis tends to infinity as $t$ tends to infinity and to zero as $t$ tends to $-\infty$. If the operator $A(t) = A$ is constant, (1) is exponentially dichotomic if and only if $A$ has no eigenvalues on the imaginary axis (see [Cop78], [DK74]).

The problem of roughness of Exponential Dichotomy with respect to small perturbations (see Section 3) has been widely studied in the scientific literature: we show that if, as has often been the case, condition (2) is assumed, then roughness is an elementary consequence of the functional analytic characterization presented in Section 1 (see [Sch99], [Min99]. [Min01]).

## 2 Characterization of Exponential Dichotomy

In order to study the functional-analytic equivalent of exponential dichotomies, we recall that given a closed unbounded, linear operator $T$ with domain $D(T)$ on a Banach space $X$, $D(T)$ becomes a Banach space when furnished with the norm $\|x\|_* = \|x\| + \|T(x)\|$ and $T : (D(t), \|\cdot\|_*) \to X$ defines a bounded linear operator. Let $S \subseteq \Re$; we will consider the Banach space

240

$C(S, X)$ consisting of all bounded, continuous functions $x : S \to X$ with the $L^\infty$-norm

$$\|x\|_{L^\infty(S,X)} = \max_{t \in S} \|x(t)\|$$

and the (unbounded) operator

$$L : D(L) = \{x \in C(\Re, X) : \frac{d}{dt}x - A(t)x \in C(\Re, X)\} \to C(\Re, X), \quad (5)$$

defined by $Lx(t) = \frac{d}{dt}x(t) - A(t)x(t)$. A straightforward argument shows that $L$ is closed. The key observation is that under the assumption of exponential growth, exponential dichotomy of the differential equation (1) is equiva- lent to the invertibility of $L$. For the sake of completeness we include the proof of this fact:

**Theorem 1.** *Under condition (2) (exponential growth), the differential equation (1) has Exponential Dichotomy if and only if the operator L introduced in (5) is invertible.*

*Proof.* For the necessity, consider the Green's function of equation (1), namely

$$G(t,s) = \begin{cases} U(t)PU^{-1}(s) \text{ if } t \geq s \\ -U(t)QU^{-1}(s) \text{ if } t \leq s \end{cases} \quad (6)$$

take $f \in C(\Re, X)$ and set $u(t) = \int_\Re G(t,s)f(s)ds$. A straightforward calculation shows that $u$ is strongly differentiable and satisfies $Lu = f$; using inequalities (3) it readily follows that

$$\|u\|_{L^\infty(\Re,X)} \leq \left(\frac{N_1}{v_1} + \frac{N_2}{v_2}\right)\|f\|_{L^\infty(\Re,X)}. \quad (7)$$

It has been established that $u \in C(\Re, X)$ and hence that $L$ is onto. The dichotomic inequalities imply that the only solution to the homogeneous equation (1) that is bounded on the real line is the trivial solution $u \equiv 0$, hence $L$ is invertible; inequality (7) shows that the inverse is bounded and that its norm is subject to the bound

241

$$\left\| L^{-1} \right\| \leq \left( \frac{N_1}{v_1} + \frac{N_2}{v_2} \right). \tag{8}$$

The sufficiency part of Theorem 1 is a consequence from Lemma 1 and Lemma 2, which we now prove:

**Lemma 1.** *If*

$$L : D(L) \to C(\Re, X) \tag{9}$$

*has a bounded inverse, then there exists positive constants $C_i$, $\alpha_i$, $i = 1,2$ depending only on the operator norm of $L^{-1}$ and $\alpha$, $\beta$ in (2) such any solution* u *of (1) for which*

$$\left\| u \right\|_{L^\infty([0,\infty)),X)} < \infty \tag{10}$$

*satisfies the estimate*

$$\left\| u(t) \right\| \leq C_1 e^{-\alpha_1(t-s)} \left\| u(s) \right\| \text{ for } t \geq s \geq 0 \tag{11}$$

*and that any solution* u *of (1) for which*

$$\left\| u \right\|_{L^\infty((-\infty,0]),X)} < \infty \tag{12}$$

*satisfies the estimate*

$$\left\| u(t) \right\| \leq C_2 e^{-\alpha_2(s-t)} \left\| u(s) \right\| \text{ for } t \leq s \leq 0. \tag{13}$$

*Proof.* Let $s \leq 0$, consider a cut-off function $\psi \in C_0^\infty(\Re)$ supported on $(-\infty, s]$ for which $\psi \equiv 1$ on $(-\infty, -s-1]$, $0 \leq \psi \leq 1$ and $\left\| \psi' \right\|_{L^\infty} \leq 2$. If $u$ is a solution of (1) then $\psi u \in D(L)$ and $L(\psi u) = \psi' u$ if in addition $u$ is bounded on the left semi-axis (i.e, satisfying the condition (12)), then by virtue of the invertibility of $L$, one has

$$\left\| u \right\|_{L^\infty((-\infty,s-1])} \leq \left\| \psi u \right\|_{L^\infty((-\infty,s-1])} \leq \left\| \psi u \right\|_{L^\infty(\Re)} \tag{14}$$

$$\leq \left\| L^{-1} \right\| \left\| \psi' u \right\|_{L^\infty(\Re)} \leq 2 \left\| L^{-1} \right\| \left\| U(t) U^{-1}(s) u(s) \right\|_{L^\infty(\Re)} \tag{15}$$

$$\leq 2\alpha e^{\beta} \|L^{-1}\| \|u(s)\|_X, \tag{16}$$

which in conjunction with the estimate

$$\sup_{t\in[s-1,s]} \|u(t)\| = \sup_{t\in[s-1,s]} \|U(t)U^{-1}(s)u(s)\| \leq \alpha e^{\beta} \|u(s)\|$$

yields the inequality

$$\|u(t)\| \leq 2\alpha\, e^{\beta} \max\{1, \|L^{-1}\|\} \|u(s)\| \text{ for } t \leq s \leq 0. \tag{17}$$

In identical fashion it can be established that any solution of (1) which is bounded on $[0,\infty)$ satisfies the inequality

$$\|u(t)\| \leq 2\alpha\, e^{\beta} \max\{1, \|L^{-1}\|\} \|u(s)\| \text{ for } t \geq s \geq 0. \tag{18}$$

Next, we consider an arbitrary interval $[a, b] \subset (-\infty, 0]$, an arbitrary solution of (1) satisfying (17) subject to $\|w(a)\| \geq 1/2$ and $\|w(b)\| \leq 1$. Then a-fortiori, for $t \in [a,b]$ and $C = 2\alpha\, e^{\beta} \max\{1, \|L^{-1}\|\}$, one has the inequality

$$\frac{1}{2C} \leq \|w(t)\| \leq C. \tag{19}$$

Take $\varepsilon > 0$ and a cut-off function $\psi$ supported on $[a, b]$, equal to 1 on $[a + \varepsilon, b - \varepsilon]$, $|\psi| \leq 1$; set

$$g(t) = \psi(t) w(t) \|w(t)\|^{-1}, \ u(t) = w(t) \int_{-\infty}^{t} \psi(s) \|w(s)\|^{-1} ds.$$

Elementary calculation show that $Lu = g$ and it follows from the assumption on $L$ that

$$\frac{1}{2c^2}(b - a - 2\varepsilon) \leq \|u\|_{L^{\infty}(\Re)} \leq \|L^{-1}\| \|g\|_{L^{\infty}(\Re)} \leq \|L^{-1}\|. \tag{20}$$

We conclude that if $u$ is any solution of (1) which is bounded on the left semi axis and $t > N > 2C^2 \|L^{-1}\|$, then the inequality

243

$$\|u(s-t)\| \leq \frac{1}{2}\|u(s)\| \tag{21}$$

holds for any $s \leq 0$. Next, let $t \leq s \leq 0$, set

$$-n-1 = \sup\{i \in Z; iN < t-s\};$$

it follows then from (21) that for some constant $0 < \varepsilon < T$ depending only on $\|L^{-1}\|$ and the exponential growth, one has

$$\|u(t)\| = \|u(t-s+s)\| \leq \frac{1}{2^n}\|u(s)\| = e^{\ln 2(-n)}\|u(s)\| \tag{22}$$

$$= e^{\ln 2(\frac{t-s}{N}+\frac{\varepsilon}{N})}\|u(s)\| \tag{23}$$

which proves (13). The proof of (11) is identical and will be omitted. Lemma (1) is proved. Next, we set

$$X_1 = \{x \in X : \exists v : v' = A(t)v, v(0) = x \text{ and } \|v\|_{L^\infty([0,\infty))} < \infty\}$$

and

$$X_2 = \{x \in X : \exists v : v' = A(t)v, v(0) = x \text{ and } \|v\|_{L^\infty((-\infty,0])} < \infty\}.$$

**Lemma 2.** *Let $X_1$ and $X_2$ be the spaces defined above. Then*

$$X = X_1 \oplus X_2$$

*algebraically and topologically. Moreover, the Differential Equation (1) is Exponentially Dichotomic with the projection on $X_1$.*

*Proof.* By virtue of (17) and (18), $X_i$, $i = 1, 2$ is a closed subspace of X. Let $\alpha \in C^\infty(\Re)$ with $\alpha \equiv 0$ on (-∞, 0], $\alpha \equiv 1$ on [1, ∞) and $\|\alpha'\|_{L^\infty(\Re)} < \infty$; for an arbitrary $x \in X$ let $u$ be the unique solution of (1) with $u(0) = x$. Set $g = \alpha' u$. Invoking again the invertibility of $L$ we can find a unique $v \in D(L)$ such that $L(v) = g$.

244

Writing

$$w = (1 - \alpha)u + v,$$

it is immediate that

$$Lw = -\alpha' u + (1 - \alpha)u' - (1 - \alpha)Au + Lv = 0;$$

in addition,

$$\|w\|_{L^\infty([0,\infty))} \le \|v\|_{L^\infty(\Re)} + \|u\|_{L^\infty([0,1])} < \infty.$$

We have proved that $w(0) \in X_1$ and since $v(0) \in X_2$, that $X=X_1+X_2$. Finally, due to the invertibility of $L$, if $z$ is bounded on the real line and $Lz=0$, then necessarily $z$ is identically 0 on $\Re$, so that the sum is direct as claimed.

To complete the proof of the Theorem we set

$$P(t) = U(t)PU^{-1}(t),\ t \in \Re$$

and show that

$$\sup_{t\in\Re}\|P(t)\| = \sup_{t\in\Re}\|U(t)PU^{-1}(t)\| < \infty. \qquad (24)$$

For $x \in X$, denote the solutions of the initial value problems $(P(t) \ne 0$ and $Q(t) \ne 0)$

$$\frac{d}{dt}u = A(t)u$$

$$u(t) = \frac{P(t)x}{\|P(t)\|}$$

and

$$\frac{d}{dt}v = A(t)v$$

$$v(t) = \frac{Q(t)x}{\|Q(t)\|}$$

245

$x = P(t)x + Q(t)x$ by $u$ and $v$ respectively and write

$$w(\xi) = \|P(t)x\|u(\xi) + \|Q(t)x\|v(\xi).$$

Hence, for $\xi \geq t$ and some $0 \leq \overline{\alpha} < \alpha$, one has

$$C\|x\|e^{\overline{\alpha}(\xi-t)} \geq \|w(\xi)\| \qquad (25)$$

$$C\|x\|e^{\overline{\alpha}(\xi-t)} \geq \|Q(t)x\|\frac{1}{C}e^{\alpha(\xi-t)} - \|P(t)\|Ce^{-\alpha(\xi-t)}$$

$$\|x\| \geq \|Q(t)x\|C_1 e^{(\alpha-\tilde{\alpha})(\xi-t)} - \|P(t)\|C_2 e^{-(\alpha+\tilde{\alpha})(\xi-t)}.$$

Recall that there exists a positive constant $c$ such that for each $\xi \geq t$, the inequalities

$$\|u(\xi)\| \leq C\|u(t)\|e^{-\alpha(\xi-t)}$$

$$\|v(\xi)\| \geq \frac{1}{c}\|v(t)\|e^{\alpha(\xi-t)}$$

hold. Choose $\xi \geq \max\left\{\frac{1}{\alpha}\ln(c)+t; t-\frac{1}{\alpha}\ln\left(\frac{\|Q(t)x\|}{2c\|P(t)x\|}\right)\right\}$.

Then, (25) yields

$$\|Q(t)x\|\frac{1}{2} \leq c\|x\|e^{\alpha(\xi-t)},$$

from which (24) follows automatically.

Lemma 2 is the final ingredient in the proof of Theorem 1.

## 3  Stability of Exponential Dichotomy

It is now easy to prove stability of Exponential Dichotomy under the assumption (2):

**Theorem 2.** *Under condition* (2)*, if* (1) *is Exponentially dichotomic, so is the perturbed equation*

$$\frac{d}{dt}x = (A(t) + B(t))x \qquad (26)$$

*for any strongly continuous operator function B : $\Re \to B(X)$ with*

$$\|B(t)\|_{L^\infty(\Re, B(X))} < \left(\frac{N_1}{v_1} + \frac{N_2}{v_2}\right)^{-1}. \qquad (27)$$

*Proof.* Since (1) is Exponentially Dichotomic, then $L$ is invertible and (9) holds; The (unbounded) operator

$$S = L - B(t) : D(L) \to C(\Re, X)$$

is invertible (see [Kat96, Sect. IV.2, Remark 2.22]). The sufficiency in Theorem 1 shows that (26) is Exponentially Dichotomic, as claimed. We underline the fact that the dichotomic constants for the perturbed equation (26) can be tracked down in the previous proof and that they depend only on the dichotomic constants of the original differential equation and the $L^\infty$ norm of $B(t)$.

*Remark 1.* Condition (2) cannot be omitted in Theorem 1, see [Cop78] for a counterexample.
 The $L^\infty$-bound (27) is optimal, see [MP08]. The roughness property of Exponential Dichotomy has been proved by methods different from the ones employed in this work, in the absence of condition (2) (see [Pop06], [MP08]).

## References

[Cop78]   W. Coppel - *Dichotomies in Stability Theory*, Springer-Verlag (1978).